\documentstyle[11pt,amssymb]{article}
\newcommand{\beql}[1]{\begin{equation}\label{#1}}
\newcommand{\eeq}{\end{equation}}
\newcommand{\comment}[1]{}

\newcommand{\eqref}[1]{{\rm (\ref{#1})}}

\newcommand{\Abs}[1]{{\left|{#1}\right|}}

\newcommand{\Norm}[1]{{\left\|{#1}\right\|}}

\newcommand{\Qed}{\ \\\mbox{$\blacksquare$}}

\newcommand{\Set}[1]{{\left\{{#1}\right\}}}

\newcommand{\ToAppear}[1]{\raisebox{15mm}[10pt][0mm]{\makebox[0mm]%
	{\makebox[\textwidth][r]{\small #1}}}}

\newcommand{\RR}{{\Bbb R}}
\newcommand{\CC}{{\Bbb C}}
\newcommand{\ZZ}{{\Bbb Z}}

\newcommand{\TT}{{\Bbb T}}

\newcommand{\one}{{\bf 1}}

\newcommand{\inner}[2]{{\langle #1, #2 \rangle}}
\newcommand{\Inner}[2]{{\left\langle #1, #2 \right\rangle}}
\newcommand{\dens}{{\rm dens\,}}

\newcommand{\supp}{{\rm supp\,}}
\newcommand{\dist}{{\rm dist\,}}
\newcommand{\vol}{{\rm vol\,}}
\newcommand{\ft}[1]{\widehat{#1}}
\newcommand{\FT}[1]{\left(#1\right)^\wedge}

\newcounter{open}
\newcounter{dfn}
\def\thedfn{\arabic{dfn}}
\newenvironment{dfn}{
  \sf
  \vskip 0.10in
  \refstepcounter{dfn}
  \noindent{\bf Definition \thedfn \ }
}{\vskip 0.10in}

\newcounter{obs}
\def\theobs{\arabic{obs}}

\newcounter{thm}
\newenvironment{thm}{
  \sf
  \vskip 0.10in
  \refstepcounter{thm}
  \noindent{\bf Theorem\ }
}{\vskip 0.10in}

\newcounter{mysec}
\def\themysec{\arabic{mysec}}
\newcommand{\mysection}[1]{
  \vskip 0.25in
  \refstepcounter{mysec}\centerline{\large\bf \S\themysec.\ {#1}}\par
  \nopagebreak
  \addcontentsline{toc}{section}{{\bf \themysec.}\ {#1}}
}

\newcounter{mysubsec}[mysec]
\newtheorem{theorem}{Theorem}

\begin{document}

\begin{center}
{\Large \bf Non-symmetric convex domains\ToAppear{
{\tt Illinois J. Math., to appear}} have no basis of exponentials}\\
\ \\
{\sc Mihail N. Kolountzakis\footnote{
Department of Mathematics, 1409 W. Green St, University of Illinois,
Urbana, IL 61801, USA. E-mail: {\tt kolount@math.uiuc.edu}}%
\footnote{Partially supported by the U.S. National Science Foundation,
under grant DMS 97-05775.}}%
\footnote{Current address: Department of Mathematics, University of Crete,
714 09 Iraklio, GREECE. E-mail: {\tt kolount@math.uch.gr}}\\
\ \\
\small December 1998; revised October 1999 \\
\end{center}

\begin{abstract}
A conjecture of Fuglede states that a bounded measurable set $\Omega\subset\RR^d$,
of measure $1$, can tile $\RR^d$
by translations if and only if the Hilbert space
$L^2(\Omega)$ has an orthonormal basis consisting of exponentials
$e_\lambda(x) = \exp 2\pi i\inner{\lambda}{x}$.
If $\Omega$ has the latter property it is called {\em spectral}.
We generalize a result of Fuglede, that a triangle in the
plane is not spectral, proving
that every non-symmetric convex domain in $\RR^d$ is not spectral.
\end{abstract}

\mysection{Introduction}\label{sec:intro}

Let $\Omega$ be a measurable subset of $\RR^d$ of measure $1$
and $\Lambda$ be a discrete subset of $\RR^d$.
We write
\begin{eqnarray*}
e_\lambda(x) &=& \exp{2\pi i \inner{\lambda}{x}},\ \ \ (x\in\RR^d),\\
E_\Lambda &=& \Set{e_\lambda:\ \lambda\in\Lambda} \subset L^2(\Omega).
\end{eqnarray*}
The inner product and norm on $L^2(\Omega)$ are
$$
\inner{f}{g}_\Omega = \int_\Omega f \overline{g},
\ \mbox{ and }\ 
\Norm{f}_\Omega^2 = \int_\Omega \Abs{f}^2.
$$

\begin{dfn}
The pair $(\Omega, \Lambda)$ is called a {\em spectral pair}
if $E_\Lambda$ is an orthonormal basis for $L^2(\Omega)$.
A set $\Omega$ will be called {\em spectral} if there is
$\Lambda\subset\RR^d$ such that 
$(\Omega, \Lambda)$ is a spectral pair.
The set $\Lambda$ is then called a {\em spectrum} of $\Omega$.
\end{dfn}

\noindent
{\bf Example:} If $Q_d = (-1/2, 1/2)^d$ is the cube of
unit volume in $\RR^d$ then
$(Q_d, \ZZ^d)$ is a spectral pair.

We write $B_R(x) = \Set{y\in\RR^d:\ \Abs{x-y}<R}$.
\begin{dfn} (Density)\\
(i) The set $\Lambda \subset \RR^d$ has {\em uniformly bounded density} if for
each $R>0$ there exists
a constant $C>0$ such that $\Lambda$ has at most $C$ elements in each ball of radius
$R$ in $\RR^d$.\\
(ii) The set $\Lambda \subset \RR^d$ has {\em density} $\rho$, and we write
$\rho = \dens \Lambda$, if we have
$$
\rho = \lim_{R\to\infty} {\Abs{\Lambda\cap B_R(x)} \over \Abs{B_R(x)}},
$$
uniformly for all $x\in\RR^d$.
\end{dfn}

We define translational tiling for complex-valued functions below.
\begin{dfn}
Let $f:\RR^d\to\CC$ be measurable and $\Lambda \subset \RR^d$ be a discrete set.
We say that {\em $f$ tiles with $\Lambda$ at level $w\in\CC$}, and sometimes
write ``$f + \Lambda = w \RR^d$'', if
\beql{tiling}
\sum_{\lambda\in\Lambda} f(x - \lambda) = w,
  \ \ \mbox{for almost every (Lebesgue) $x\in\RR^d$},
\eeq
with the sum above converging absolutely a.e.
If $\Omega \subset \RR^d$ is measurable we say that {\em $\Omega + \Lambda$ is a tiling}
when $\one_\Omega + \Lambda = w\RR^d$, for some $w$.
If $w$ is not mentioned it is understood to be equal to $1$.
\end{dfn}

\noindent
{\bf Remarks}\\
1. If $f\in L^1(\RR^d)$ and $\Lambda$ has uniformly
bounded density one can easily show
(see \cite{KL96} for the proof in one dimension,
which works in higher dimension as well) that
the sum in \eqref{tiling} converges absolutely a.e.\ and
defines a locally integrable function of $x$.\\
2. In the very common case when $f\in L^1(\RR^d)$ and
$\int_{\RR^d} f\neq 0$ the condition that
$\Lambda$ has uniformly bounded density follows easily
from \eqref{tiling} and need not be postulated a priori.\\
3. It is easy to see that if $f\in L^1(\RR^d)$, $\int_{\RR^d}f \neq 0$
and $f+\Lambda$ is a tiling then $\Lambda$ has a density
and the level of the tiling $w$ is given by
$$
w = \int_{\RR^d}f \cdot \dens \Lambda.
$$

From now on we restrict ourselves to tiling with functions in $L^1$ and sets of finite measure.

\noindent
{\bf Example:} $Q_d + \ZZ^d$ is a tiling.

The following conjecture is still unresolved.

\noindent
{\bf Conjecture:} (Fuglede \cite{F74}) If $\Omega \subset \RR^d$
is bounded and has Lebesgue measure $1$ then $L^2(\Omega)$ has
an orthonormal basis of exponentials if and only if
there exists $\Lambda \subset \RR^d$ such that $\Omega + \Lambda = \RR^d$ is a tiling.

\noindent{\bf Remark:}
It is not hard to show \cite{F74} that $L^2(\Omega)$
has a basis $\Lambda$ which
is a {\em lattice} (i.e., $\Lambda = A \ZZ^d$, where $A$
is a non-singular $d\times d$ matrix)
if and only if $\Omega+\Lambda^*$ is a tiling.
Here
$$
\Lambda^* = \Set{\mu\in\RR^d:\ \inner{\mu}{\lambda} \in \ZZ,\ \forall \lambda\in\Lambda}
$$
is the {\em dual lattice} of $\Lambda$ (we have $\Lambda^* = A^{-\top}\ZZ^d$).

\comment{
Recently Lagarias, Reeds and Wang \cite{LRW98} and Iosevich and Pedersen \cite{IP98}
have shown that $E_\Lambda$ is an orthnormal basis for $L^2(Q_d)$ if and only if
$Q_d + \Lambda$ is a tiling. The interest of this result lies in the fact that $\Lambda$
may well not be a lattice.
}

Fuglede \cite{F74} showed that the disk
and the triangle in $\RR^2$ are not spectral domains.
\comment{
He also showed that $L^2$ of the triangle contains
an infinite orthogonal family of exponentials.
}

In this note we prove the following generalization of Fuglede's triangle result.
\begin{theorem}\label{th:nonsym}
Let $\Omega$ have measure $1$ and be a convex, non-symmetric, bounded open set in $\RR^d$.
Then $\Omega$ is not spectral.
\end{theorem}

The set $\Omega$ is called {\em symmetric}
with respect to $0$ if $y \in \Omega$ implies
$-y \in \Omega$, and symmetric with respect
to $x_0\in\RR^d$ if $y\in\Omega$ implies
that $2x_0-y\in\Omega$.
It is called {\em non-symmetric} if it is not
symmetric with respect to any $x_0\in\RR^d$.
For example, in any dimension a simplex is non-symmetric.

It is known \cite{V54,M80} that every convex body that tiles $\RR^d$
by translation is a centrally symmetric polytope and that each such body
also admits a lattice tiling and, therefore
(see the remark after Fuglede's conjecture above), its
$L^2$ admits a lattice spectrum. Given Theorem \ref{th:nonsym},
to prove Fuglede's conjecture restricted to convex domains,
one still has to prove that any symmetric convex body that is not
a tile admits no orthonormal basis of exponentials for its $L^2$.

In \S\ref{sec:fourier} we derive some necessary
and some sufficient conditions for
$f+\Lambda$ to be a tiling. 
These conditions roughly state that tiling is equivalent to
a certain tempered distribution, associated with $\Lambda$
being ``supported'' on the zero set of $\ft{f}$ plus the origin.
Similar conditions had been derived in \cite{KL96} but
here we have to work with less smoothness for $\ft{f}$.
To compensate for the lack of smoothness we work with compactly
supported $\ft{f}$ and nonnegative $f$ and $\ft{f}$,
conditions which are fulfilled for our problem.

In \S\ref{sec:main} we restate the property
that $\Omega$ is spectral as a tiling problem
for $\Abs{\ft{\one_\Omega}}^2$ and use the conditions
derived in \S\ref{sec:fourier} to prove Theorem \ref{th:nonsym}.
What makes the proof work is that when $\Omega$ is a non-symmetric
convex set the set $\Omega-\Omega$ has volume
strictly larger than $2^d \vol\Omega$.

\mysection{Fourier-analytic conditions for tiling}\label{sec:fourier}

Our method relies on a Fourier-analytic characterization of translational tiling,
which is a variation of the one used in \cite{KL96}.
We define the (generally unbounded) measure
$$
\delta_\Lambda = \sum_{\lambda\in\Lambda} \delta_\lambda,
$$
where $\delta_\lambda$ represents a unit mass at $\lambda\in\RR^d$.
If $\Lambda$ has uniformly bounded density then $\delta_\Lambda$
is a tempered distribution
(see for example \cite{R73}) and therefore its Fourier Transform
$\ft{\delta_\Lambda}$ is defined and is itself a tempered distribution.

The action of a tempered distribution (see \cite{R73})
$\alpha$ on a Schwartz function $\phi$
is denoted by $\alpha(\phi)$. The Fourier Transform of $\alpha$ is
defined by the equation
$$
\ft\alpha(\phi) = \alpha(\ft\phi).
$$
The support $\supp\alpha$ is the smallest closed set $F$ such that for any
smooth $\phi$ of compact support contained in the open set $F^c$ we have
$\alpha(\phi)=0$.
\begin{theorem}\label{th:tiling-implies-supp}
Suppose that $f\ge 0$ is not identically $0$, that $f \in L^1(\RR^d)$,
$\widehat{f}\ge 0$ has compact support and $\Lambda\subset\RR^d$.
If $f+\Lambda$ is a tiling then
\beql{sp-cond-1}
\supp \ft{\delta_\Lambda} \subseteq \Set{x\in\RR^d:\ \ft{f}(x) = 0} \cup \Set{0}.
\eeq
\end{theorem}
{\bf Proof of Theorem \ref{th:tiling-implies-supp}.}
Assume that $f+\Lambda = w\RR^d$ and let
$$
K = \Set{\widehat{f} = 0} \cup \Set{0}.
$$
We have to show that
$$
\widehat{\delta_\Lambda}(\phi) = 0,
\ \ \forall \phi\in C_c^\infty(K^c).
$$
Since $\ft{\delta_\Lambda}(\phi) = \delta_\Lambda(\ft\phi)$
this is equivalent to $\sum_{\lambda\in\Lambda}\widehat\phi(\lambda) = 0$, for
each such $\phi$.
Notice that $h = \phi / \widehat{f}$ is a continuous function, but
not necessarily smooth. We shall need that $\widehat h \in L^1$.
This is a consequence of a well-known theorem of Wiener \cite[Ch.\ 11]{R73}.
We denote by $\TT^d = \RR^d/\ZZ^d$ the $d$-dimensional torus.
\begin{thm} (Wiener)\\
If $g \in C(\TT^d)$ has an absolutely convergent Fourier series
$$
g(x) = \sum_{n\in\ZZ^d} \widehat g(n) e^{2\pi i \Inner{n}{x}},
\ \ \ \widehat g \in \ell^1(\ZZ^d),
$$
and if $g$ does not vanish anywhere on $\TT^d$ then
$1/g$ also has an absolutely convergent Fourier series.
\end{thm}
Assume that
$$
\supp\phi,\ \supp\widehat{f} \subseteq \left(-{L\over2},{L\over2}\right)^d.
$$
Define the function $F$ to be:\\
(i) periodic in $\RR^d$ with period lattice $(L\ZZ)^d$,\\
(ii) to agree with $\widehat{f}$ on $\supp\phi$,\\
(iii) to be non-zero everywhere and,\\
(iv) to have $\widehat{F} \in \ell^1(\ZZ^d)$, i.e.,
$$
\widehat F = \sum_{n\in\ZZ^d} \widehat F(n) \delta_{L^{-1} n},
$$
is a finite measure in $\RR^d$.

One way to define such an $F$ is as follows.
First, define the $(L\ZZ)^d$-periodic function $g\ge 0$
to be $\ft{f}$ periodically extended.
The Fourier coefficients of $g$ are $\ft{g}(n) = L^{-d} f(-n/L) \ge 0$.
Since $g, \ft g \ge 0$ and $g$ is continuous at $0$ it is
easy to prove that $\sum_{n\in\ZZ^d} \ft g (n) = g(0)$, and therefore
that $g$ has an absolutely convergent Fourier series.

Let $\epsilon$ be small enough to guarantee that $\ft{f}$
(and hence $g$) does not vanish on $(\supp\phi) + B_\epsilon(0)$.
Let $k$ be a smooth $(L\ZZ)^d$-periodic function which is equal to $1$
on $(\supp\phi)+(L\ZZ^d)$ and equal to $0$ off
$(\supp\phi + B_\epsilon(0))+(L\ZZ^d)$, and
satisfies $0\le k \le 1$ everywhere.
Finally, define
$$
F = k g + (1-k).
$$
Since both $k$ and $g$ have absolutely summable Fourier series and this
property is preserved under both sums and products, it follows that $F$ also
has an absolutely summable Fourier series. And by the nonnegativity of $g$ we
get that $F$ is never $0$, since $k=0$ on $Z(\ft f)+(L\ZZ^d)$.

By Wiener's theorem, $\widehat{F^{-1}} \in \ell^1(\ZZ^d)$, i.e., $\ft{F^{-
1}}$ is a finite
measure on $\RR^d$.
We now have that
$$
\FT{{\phi \over \widehat f}} = \widehat{\phi F^{-1}} =
  \widehat\phi * \widehat{F^{-1}} \in L^1(\RR^d).
$$
This justifies the interchange of the summation and integration below:
\begin{eqnarray*}
\sum_{\lambda\in\Lambda} \widehat\phi(\lambda) 
 &=& \sum_{\lambda\in\Lambda} \FT{{\phi\over\widehat{f}} \widehat{f}}
(\lambda) \\
 &=& \sum_{\lambda\in\Lambda} \FT{{\phi\over\widehat f}} *
\widehat{\widehat{f}}~(\lambda) \\
 &=& \sum_{\lambda\in\Lambda} \int_{\RR^d} \FT{{\phi\over\widehat f}}(y) f(y-
\lambda) ~dy \\
 &=& \int_{\RR^d} \FT{{\phi\over\widehat f}}(y)
  \sum_{\lambda\in\Lambda} f(y-\lambda) ~dy\\
 &=& w \int_{\RR^d} \FT{{\phi\over\widehat f}}(y)~dy\\
 &=& w{\phi\over\ft{f}}(0)\\
 &=& 0,
\end{eqnarray*}
as we had to show.
\Qed

For a set $A \subseteq \RR^d$ and $\delta > 0$ we write
$$
A_\delta = \Set{x\in\RR^d:\ \dist(x, A)<\delta}.
$$

We shall need the following partial converse to Theorem
\ref{th:tiling-implies-supp}.
\begin{theorem}\label{th:disjoint-supports}
Suppose that $f\in L^1(\RR^d)$,
and that $\Lambda \subset \RR^d$ has uniformly bounded density.
Suppose also that $O \subset \RR^d$ is open and
\beql{dijoint-supports}
\supp \ft{\delta_\Lambda} \setminus \Set{0} \ \subseteq\ O
\ \mbox{ and }\ O_\delta \subseteq\  \Set{\ft{f} = 0},
\eeq
for some $\delta>0$.
Then $f+\Lambda$ is a tiling at level $\ft{f}(0)\cdot\ft{\delta_\Lambda}(\Set{0})$.
\end{theorem}
{\bf Proof.}
Let $\psi:\RR^d\to\RR$ be smooth, have support in $B_1(0)$ and $\ft\psi(0)=1$ and for
$\epsilon>0$ define the approximate identity
$\psi_\epsilon(x) = \epsilon^{-d}\psi(x/\epsilon)$.
Let
$$
f_\epsilon = \ft{\psi_\epsilon} f,
$$
which has rapid decay.

First we show that $(\int f_\epsilon)^{-1}f_\epsilon + \Lambda$
is a tiling.
That is, we show that the convolution $f_\epsilon * \delta_\Lambda$
is a constant.
Let $\phi$ be any Schwartz function. Then
$$
f_\epsilon * \delta_\Lambda (\phi) =
\ft{f_\epsilon}\ft{\delta_\Lambda} (\ft\phi(-x)) =
\ft{\delta_\Lambda} (\ft\phi(-x) \ft{f_\epsilon}).
$$
The function $\ft\phi(-x) \ft{f_\epsilon}$ is a Schwartz
function whose support intersects $\supp\ft{\delta_\Lambda}$ only at $0$,
since, for small enough $\epsilon>0$,
$$
\supp \ft\phi \ft{f_\epsilon} \subseteq
\supp \ft{f_\epsilon} \subseteq
(\supp\ft{f})_\epsilon \subseteq O^c.
$$
Hence, for each Schwartz function $\phi$
$$
f_\epsilon * \delta_\Lambda (\phi) = 
	\ft\phi(0) \ft{f_\epsilon}(0) \ft{\delta_\Lambda}(\Set{0}),
$$
which implies
$$
f_\epsilon * \delta_\Lambda (x) = \ft{f_\epsilon}(0) \ft{\delta_\Lambda}(\Set{0}),
\ \ \mbox{a.e.($x$)}.
$$
We also have that $\sum_{\lambda\in\Lambda} \Abs{f(x-\lambda)}$ is finite
a.e.\ (see Remark 1 following the definition of tiling), hence,
for almost every $x\in\RR^d$
$$
\sum_{\lambda\in\Lambda}\Abs{f(x-\lambda) - f_\epsilon(x-\lambda)} =
\sum_{\lambda\in\Lambda}\Abs{f(x-\lambda)}\cdot\Abs{1-\ft{\psi_\epsilon}(x-\lambda)},
$$
which tends to $0$ as $\epsilon\to 0$.
This proves
$$
\sum_{\lambda\in\Lambda} f(x-\lambda) =
	\ft{f}(0)\cdot\ft{\delta_\Lambda}(\Set{0}),\ \ \mbox{a.e.($x$)}.
$$
\Qed

\mysection{Proof of the main result}\label{sec:main}

We now make some remarks that relate the property of $E_\Lambda$
being a basis for $L^2(\Omega)$ to a certain function tiling $\RR^d$ with $\Lambda$.

Assume that $\Omega$ is a bounded open set of measure $1$.
Notice first that
$$
\inner{e_\lambda}{e_x}_\Omega = \ft{\one_\Omega}(x-\lambda).
$$
The set $E_\Lambda$ is an orthonormal basis for $L^2(\Omega)$ if and only if for each
$f\in L^2(\Omega)$
$$
\Norm{f}_\Omega^2 = \sum_{\lambda\in\Lambda} \Abs{\inner{e_\lambda}{f}_\Omega}^2,
$$
and, by the completeness of the exponentials in $L^2$ of a large cube containing $\Omega$,
it is necessary and sufficient that
\beql{tiling-condition-for-spectrum}
\sum_{\lambda\in\Lambda} \Abs{\ft{\one_\Omega}(x-\lambda)}^2 = 1,
\eeq
for each $x\in\RR^d$.
In other words a necessary and sufficient condition for $(\Omega,\Lambda)$ to be a spectral
pair is that
$\Abs{\ft{\one_\Omega}}^2 + \Lambda$ is a tiling at level $1$.
Notice also that $\Abs{\ft{\one_\Omega}}^2$ is the Fourier Transform of
$\one_\Omega * \widetilde{\one_\Omega}$ which has support equal to the set
$\overline{\Omega - \Omega}$.
We use the notation $\widetilde{f}(x) = \overline{f(-x)}$.

\noindent
{\bf Proof of Theorem \ref{th:nonsym}:}
Write $K = \Omega - \Omega$, which is a symmetric, open convex set.
Assume that $(\Omega, \Lambda)$ is a spectral pair.
We can clearly assume that $0 \in \Lambda$.
It follows that $\Abs{\ft{\one_\Omega}}^2 + \Lambda$ is a tiling and hence that 
$\Lambda$ has uniformly bounded density, has density equal to $1$ and
$\ft{\delta_\Lambda}(\Set{0}) = 1.$

By Theorem \ref{th:tiling-implies-supp} (with
$f = \Abs{\ft{\one_\Omega}}^2,\ \ \ft{f} =
	\one_\Omega * \widetilde{\one_\Omega}(-x)$)
it follows that
$$
\supp \ft{\delta_\Lambda} \subseteq \Set{0} \cup K^c.
$$
Let $H = K/2$ and write
$$
f(x) = \one_H*\widetilde{\one_H}(x) = \int_{\RR^d} \one_H(y) \one_H(y-x)~dy.
$$
The function $f$ is supported in $\overline{K}$ and has nonnegative Fourier Transform
$$
\ft{f} = \Abs{\ft{\one_H}}^2.
$$
We have
$$
\int_{\RR^d}\ft{f} = f(0) = \vol H
$$
and
$$
\ft{f}(0) = \int_{\RR^d} f = (\vol H)^2.
$$
By the Brunn-Minkowski inequality (see for example \cite[Ch.\ 3]{G94}),
for any convex body $\Omega$,
$$
\vol {1\over2}(\Omega-\Omega) \ge \vol \Omega,
$$
with equality only in the case of symmetric $\Omega$.
Since $\Omega$ has been assumed to be non-symmetric it follows that
$$
\vol H > 1.
$$
For
$$
1>\rho>\left({1\over \vol H}\right)^{1/d}
$$
consider
$$
g(x) = f(x/\rho)
$$
which is supported properly inside $K$, and has
$$
g(0) = f(0) = \vol H,\ \ \int_{\RR^d} g = \rho^d \int_{\RR^d} f = \rho^d (\vol H)^2.
$$
Since $\supp g$ is properly contained in $K$
Theorem \ref{th:disjoint-supports} implies that $\ft{g} + \Lambda$ is a tiling
at level $\int\ft{g} \cdot \dens\Lambda = \int\ft{g} = g(0) = \vol H$.
However, the value of $\ft{g}$ at $0$ is $\int g = \rho^d (\vol H)^2 > \vol H$,
and, since $\ft{g} \ge 0$ and $\ft{g}$ is continuous, this is a contradiction.
\Qed

\mysection{Bibliography}

\end{document}